\documentclass[11pt]{article}

\usepackage{amsmath,graphicx,dsfont}
\usepackage{array}
\usepackage{epsfig}
\usepackage{graphicx}
\usepackage{psfrag}
\usepackage{amsmath}
\usepackage{amssymb,subfigure}
\usepackage{mathrsfs}
\usepackage{color}
\usepackage{pstricks,pst-node,pst-text,pst-3d}
\usepackage{pst-plot}

\def\R{\mathbb R}

\def\bu{\overline{u}}
\def\bv{\overline{v}}
\def\uu{\underline{u}}
\def\uv{\underline{v}}
\def\uut{\underline{u}_{\theta}}
\def\uuto{\underline{u}_{\theta,0}}

\def\tom{\overline{t}_{\omega}}
\def\lp {\left( }
\def\rp {\right) }

\def\ie{\emph{i.e. }}

\def\epsilon{\varepsilon}
\def\ds{\displaystyle}
\def\ph{\phantom{1}}

\newcommand{\be}{\begin{equation}}
\newcommand{\ee}{\end{equation}}
\newcommand{\baa}{\begin{array}}
\newcommand{\eaa}{\end{array}}
\newcommand{\ba}{\begin{eqnarray}}
\newcommand{\ea}{\end{eqnarray}}

\newcommand\C[1]{\mathcal{C}^{#1}}
\newcommand{\carre}{\hfill$\Box$\par\addvspace{4mm}}

\newtheorem{theo}{Theorem}

\newtheorem{hyp}{Hypothesis}

\title{Accelerating solutions in integro-differential equations}

\author{Jimmy Garnier$^{\hbox{ \small{a,b}}}$ \\
\\
\footnotesize{$^{\hbox{a }}$UR 546 Biostatistique et Processus Spatiaux, INRA, F-84000 Avignon, France}\\
\footnotesize{$^{\hbox{b }}$Aix-Marseille Universit\'e, LATP, Facult\'e des Sciences et Techniques}\\
\footnotesize{Avenue Escadrille Normandie-Niemen, F-13397 Marseille Cedex 20, France}\\
}

\date{}

\begin{document}

\maketitle

\begin{abstract}

In this paper, we study the spreading properties of the solutions of an integro-differential equation of the form $u_t=J\ast u-u+f(u).$ We focus on equations with slowly decaying dispersal kernels $J(x)$ which 
correspond to models of population dynamics with long-distance dispersal events. We prove that for kernels $J$ which decrease to $0$ slower than any exponentially decaying function, the 
level sets of the solution $u$ propagate with an infinite asymptotic speed. Moreover, we obtain lower and upper bounds for the position of any level set of $u.$ These bounds 
 allow us to estimate how the solution accelerates, depending on the kernel $J$: the slower the kernel decays, the faster the level sets propagate.
Our results are in sharp contrast with most results on this type of equation, where the dispersal kernels are generally assumed to  decrease exponentially fast, leading 
to finite propagation speeds.  
\end{abstract}


\noindent{\it AMS\/}:
47G20, 45G10, 35B40.


\section{Introduction and main assumptions}

In this paper we study the large-time behavior of the solutions of  integro-differential equations with slowly decaying dispersal kernels. Namely, we consider the Cauchy problem:
\be
\left\{
           \begin{array}{l}
            u_t=J\ast u- u +f(u), \ t>0,\ x\in\R\\
            u(0,x)=u_0(x),\ x\in\R
           \end{array}
\right.
\label{eq:1}
\ee
where $J(x)$ is the dispersal kernel and
\[
(J\ast u)(t,x)=\int_{\R}{J(x-y)u(t,y)dy}.        
\]
We assume that the nonlinearity $f$ is monostable and that the initial condition $u_0$ is compactly supported. 

The equation \eqref{eq:1} arises in population dynamics~\cite{Fif79b,MedKot03} where the unknown quantity $u$ typically stands for a population density. One of the most interesting features of this model,
compared to reaction-diffusion equations, is that it can take rare long-distance dispersal events into account. Therefore, equation \eqref{eq:1} and other closely related equations have been used to explain
some rapid propagation phenomena that could hardly be explained with reaction-diffusion models, at least with compactly supported initial conditions. A classical example is Reid's paradox of rapid plant 
migration~\cite{Cla98,ClaFas98,Ske51} which is usually explained using integro-differential equations with slowly decaying kernels or with reaction-diffusion equations with slowly decaying -- and therefore 
noncompact -- initial conditions~\cite{RoqHamFayFad10}. As we shall see in this paper, the use of slowly decaying dispersal kernels is the key assumption that leads to qualitative behavior of the 
solution of~\eqref{eq:1} very different from what is expected with reaction-diffusion equations.

Let us make our assumptions more precise. We assume that the initial condition $u_0 : \R \rightarrow [0, 1]$
is continuous, compactly supported and not identically equal to $0$. 

The reaction term $f : [0, 1]\rightarrow \R$ is of class $\C{1}$ and satisfies:
\be
 f(0)=f(1)=0, \ph f(s)>0 \ph \hbox{for all}\ph s\in(0,1), \hbox{ and } f'(0)>0.
\label{hyp:f}
\ee
A particular class of such reaction term is that of Fisher-KPP type~\cite{Fis37,kpp}. For this class, the growth rate $f(s)/s$ is maximal at $s = 0.$
Furthermore, we assume that there exist $\delta > 0,$ $s_0 \in (0, 1)$ and $M \geq 0$ such that
\begin{equation}\label{hyp:f1}
 f (s) \geq f'(0) s - M s^{1+\delta} \hbox{ for all } s \in [0, s_0 ].
\end{equation}
This last assumption is readily satisfied if $f$ is of class $C^{1,\delta}.$

We assume that the kernel $J : \R \rightarrow \R$ is a nonnegative even function of mass one and with finite first moment:
\be
 J\in\C{0}(\R),\ph J>0, \ph J(x)=J(-x), \ph \int_{\R}{J(x)dx}=1 \hbox{ and } \int_{\R}{|x| J(x)dx}<\infty.
\label{hyp:J1}
\ee
Furthermore, we assume that $J(x)$ is decreasing for all $x\geq0,$ $J$ is a $\C1$ function for large $x$ and 
\begin{equation}\label{hyp:J2}
 J'(x)=o(J(x))  \hbox{ as } |x|\to+\infty.
\end{equation}
This last condition implies
that $J$ decays more slowly than any exponentially decaying functions as $|x| \rightarrow \infty$, in the sense that
\be
 \forall \eta>0,\ph \exists \ph x_{\eta}\in\R,\ph J(x)\geq e^{-\eta x}\ph \hbox{in}\ph [x_{\eta},\infty),
\label{hyp:J3}
\ee
or, equivalently, $J(x)e^{\eta |x|} \rightarrow\infty$ as $|x|\rightarrow\infty$ for all $\epsilon > 0$. We shall refer to functions $J$ satisfying the above 
assumptions~\eqref{hyp:J1},~\eqref{hyp:J2} as \textit{exponentially unbounded kernels}.

The assumption~\eqref{hyp:J2} is in contrast with the large mathematical literature on
integro-differential equations~\cite{Aro77,CovDup07,Die79,Thi79,Wei82,Wei02} as well as integro-difference equations~\cite{Lut07,Lut08}, where the dispersal kernels $J$ are generally assumed to be exponentially 
bounded as $|x|\rightarrow\infty,$ \ie:
\be
\exists \, \eta>0 \hbox{ such that } \int_{\R}J(x) e^{\eta |x|}<\infty.
\ee
In this ``exponentially bounded case", it follows from the results in~\cite{Wei82} that, under our assumptions on $u_0$ and $f$, the solution of~\eqref{eq:1} admits a \textit{finite} spreading speed $c^*.$ Thus,
 for any $c_1,$ $c_2$ with
$0 < c_1 < c^* < c_2 < \infty $ the solution $u$ to~\eqref{eq:1} tends to zero uniformly in the region $|x|\geq c_2 t,$ whereas it is bounded away from zero uniformly
in the region $|x|\leq c_1t$ for $t$ large enough.  Thus, the spreading properties of the solution of~\eqref{eq:1} when $J$ is exponentially bounded are quite similar to that of the solution of the
 reaction-diffusion equation $u_t=u_{xx}+f(u)$ with $u(0,\cdot)=u_0$~\cite{AroWei75,AroWei78,Fis37,kpp}.
The existence of such a finite spreading speed is also true for other integro-differential equations with exponentially bounded dispersal kernels~\cite{Aro77,Die79,Thi79,Wei82}.

Let us come back to problem~\eqref{eq:1} with  an exponentially unbounded kernel $J.$ In this case, it is known that equation \eqref{eq:1} does not admit any traveling wave solution with constant speed
and constant (or periodic) profile~\cite{Yag09}. Moreover, numerical results and formal analytic computations carried out for linear integro-difference equations~\cite{KotLew96} and linear
integro-differential equations~\cite{MedKot03} indicate that exponentially unbounded dispersal kernels
 lead to accelerating propagation phenomena and infinite spreading speeds.
In this article, we prove rigorously such results for the solution $u$ of~\eqref{eq:1} when the kernel $J$ is exponentially unbounded \ie $J$ satisfies~\eqref{hyp:J1} and \eqref{hyp:J2}.

Our approach is inspired from~\cite{HamRoqP2}, where it was shown for a reaction-diffusion equation $u_t=u_{xx}+f(u)$ that exponentially unbounded initial conditions lead to solutions which accelerate and have
infinite spreading speed.
Here, we get comparable results starting from compactly supported initial data and with exponentially
 unbounded dispersal kernels. However, the interpretation of our results as well as their proofs are very different from those in~\cite{HamRoqP2,RoqHamFayFad10}.
These differences are mostly due to the nonlocal nature of the operator $u\mapsto \,J\ast u - u,$ and to its lake of regularization properties.

\section{Main results}\label{sec:main_results}
 Before stating our main results, we recall that from the maximum
principle~\cite{Wei82,Yag09} and from the assumptions on $u_0,$ the solution $u$ of~\eqref{eq:1} satisfies
\[
0 < u(t, x) < 1 \ph\hbox{for all}\ph t > 0\ph \hbox{and}\ph x \in \R.
\]
For any $\lambda\in (0, 1)$ and $t \geq 0$ we denote by
\[
 E_\lambda(t)=\{x\in\R, \ph u(t,x)=\lambda\},
\]
the level set of $u$ of value $\lambda$ at time $t$. For any subset $A \subset (0, J(0))$, we set
\[
  J^{-1}\{A\}=\{x\in\R, \ph J(x)\in A\},
\]
the inverse image of $A$ by $J$.

Our first result says that the level sets $E_\lambda(t)$ of all level values $\lambda \in (0, 1)$ (namely, the time-dependent sets of real numbers $x$ such that $u(t, x) = \lambda$) move infinitely fast
as $t \to\infty.$
\begin{theo}\label{theo:infinite_speed}
Let $u$ be the solution of~\eqref{eq:1} with a continuous and compactly supported initial condition $u_0: \, \R \to [0,1]$ $(u_0\not\equiv0).$ Assume that $J$ is an exponentially unbounded kernel 
satisfying~\eqref{hyp:J1} and~\eqref{hyp:J2}.
Then, 
\begin{equation}\label{eq:E_lin}
 \forall c\geq0,\ \min_{|x|\leq ct}{u(t,x)}\to0 \hbox{ as } t\to\infty
\end{equation}

and for any given $\lambda\in(0, 1)$, there is a real number $t_\lambda \geq 0$ such that $E_\lambda(t)$ is non-empty for all $t \geq t_\lambda$ , and
\be\label{eq:E}
 \lim_{t\rightarrow+\infty}{\frac{\min \{ E_\lambda(t)\cap[0,+\infty)\}}{t}}=\lim_{t\rightarrow+\infty}{\frac{-\max \{ E_\lambda(t)\cap(-\infty,0]\}}{t}}=+\infty.
\ee
\end{theo}

Our next result gives a ``lower bound" for the level sets $E_\lambda(t)$ in terms of the behavior of $J$ at $\infty.$
\begin{theo}\label{theo:speed_estimate_inf}
Under the same asumptions as in Theo.\ref{theo:infinite_speed}, for any $\lambda\in(0,1)$ and $\epsilon\in(0,f'(0))$ there exists $T_{\lambda,\epsilon}\geq t_\lambda$ such that
\begin{equation}\label{ineq:E_inf}
    \forall \, t\geq T_{\lambda,\epsilon}, \ E_\lambda(t)\subset J^{-1}\left\{\left( 0 , e^{-(f'(0)-\epsilon)t} \right]\right\}.
\end{equation}
\end{theo}

In our next result, we will either assume:
\begin{hyp}\label{hyp:J_nonexp}
An exponentially unbounded kernel $J$ satisfies Hypothesis~\ref{hyp:J_nonexp} if and only if there exists $\sigma>0$ such that $|J'(x)/J(x)|$ is nonincreasing for all $x\geq\sigma$
and there exists $\epsilon_0\in(0,1)$ such that
\begin{equation}\label{eq:J_nonexp}
 \int_{\R}{J(z)^{\epsilon_0}dz}<\infty.
\end{equation}
\end{hyp}
or
\begin{hyp}\label{hyp:J_geo}
An exponentially unbounded kernel $J$ satisfies Hypothesis~\ref{hyp:J_geo} if and only if
\begin{equation}\label{eq:J_geo}
\ds\left|\frac{J'(x)}{J(x)}\right|=O\lp\frac{1}{|x|}\rp \hbox{ as }|x|\to\infty.
\end{equation}
\end{hyp}

Under these additional assumptions on the kernel $J,$ we are able to establish an ``upper bound" for the level sets $E_\lambda(t).$
\begin{theo}\label{theo:speed_estimate_sup}
Let $u$ be the solution of~\eqref{eq:1} with a continuous and compactly supported initial condition $u_0: \, \R \to [0,1]$ $(u_0\not\equiv0).$ Assume that $J$ satisfies either Hyp.~\ref{hyp:J_nonexp} or 
Hyp.~\ref{hyp:J_geo}.
Then, there exists $\rho>f'(0)$ such that for any $\lambda\in(0,1)$ there is $T_{\lambda}\geq t_\lambda$ such that
\begin{equation}\label{ineq:E_sup}
    \forall t\geq T_{\lambda}, \ E_\lambda(t)\subset J^{-1}\left\{\left[ e^{-\rho t}, J(0) \right]\right\}.
\end{equation}
\end{theo}

Theorem~\ref{theo:speed_estimate_inf} together with Theorem~\ref{theo:speed_estimate_sup} provide an estimation of the position of the level sets $E_\lambda(t)$ for large time $t.$ 
In particular the inclusions~\eqref{ineq:E_inf} and~\eqref{ineq:E_sup} mean that, for any $\lambda\in(0, 1)$ and any element $x_\lambda(t) \in E_\lambda(t),$ we have
\begin{equation}
 \min{\lp J^{-1}\lp  e^{-(f'(0)-\epsilon) t}\rp\cap[0,+\infty)\rp} \leq\ |x_\lambda(t)|\ \leq \max{\lp J^{-1}\lp e^{-\rho t}\rp\cap[0,+\infty)\rp},
\end{equation}
for large $t.$

\section{Case studies}
Let us apply the results of Sec.~\ref{sec:main_results} to several examples of exponentially unbounded kernels:
\begin{itemize}
\item Functions $J$ which are logarithmically sublinear as $|x|\rightarrow\infty$, that is
\be\label{exp:NEpLog}
 \ds J(x)= Ce^{-\alpha |x|/\ln(|x|)} \hbox{ for large } |x|,
\ee
with $\alpha>0$, $C > 0$;

\item Functions $J$ which are logarithmically power-like and sublinear as $|x|\rightarrow\infty$, that is
\be\label{exp:NEp}
 \ds J(x)= Ce^{-\beta|x|^{\alpha}} \hbox{ for large } |x|,
\ee
with $\alpha\in(0,1)$, $\beta,$ $C > 0$;

\item Functions $J$ which decay algebraically as $|x|\rightarrow\infty$, that is
\be\label{exp:Geo}
 \ds J(x)= C|x|^{-\alpha} \hbox{ for large } |x|,
\ee
with $\alpha>2$, $C > 0$.
.
\end{itemize}

First, if $J$ satisfies~\eqref{hyp:J1} and~\eqref{exp:NEpLog}
then $J$ satisfies Hyp.~\ref{hyp:J_nonexp} (but not Hyp.~\ref{hyp:J_geo}). Theorem~\ref{theo:speed_estimate_inf} and~\ref{theo:speed_estimate_sup} then imply that for any level value $\lambda\in(0,1)$ 
and any $\epsilon>0,$ it exists $\tilde \rho > f'(0)$ such that every element $x_\lambda(t)$ in the level set $E_\lambda(t)$ satisfies:
\begin{equation}\label{ineq:xl_NEpLog}
 \frac{f'(0)-\epsilon}{\alpha}\ t\ln{(t)}\ \leq\  |x_\lambda(t)| \ \leq \ \frac{\tilde \rho}{\alpha}t\ln(t) \ \hbox{ for large } t.
\end{equation}

\begin{figure}[ht]
\begin{center}
\includegraphics[width=8cm]{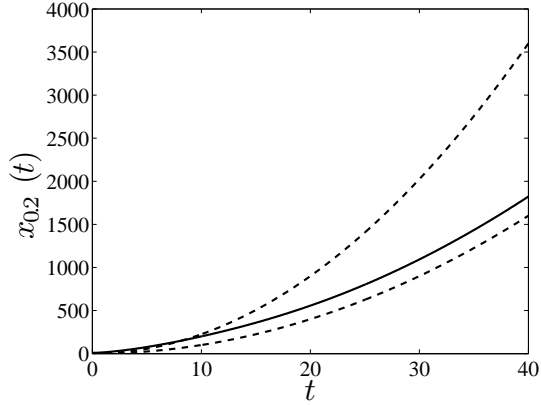}
\end{center}
\caption{Plain line: position  $x_{0.2}(t)$ of the level set $E_{0.2}(t)$ of the solution of~\eqref{eq:1} with $f(u)=u(1-u),$ $u_0(x)=\max((1-(x/10)^2),0)$ and the exponentially unbounded kernel 
$J(x)=(1/4) e^{-\sqrt{|x|}}$. Observe that $x_{0.2}(t)$ remains bounded by  $t\mapsto J^{-1}(e^{-f'(0)t})=(t-\ln(4))^2$ and $J^{-1}(e^{-(f'(0)+1/2) t})=( 3t/2-\ln(4))^{2}$ (dashed lines) for large $t.$ }
\label{fig:xl_NEp}
\end{figure}

Now, if $J$ satisfies~\eqref{hyp:J1} and~\eqref{exp:NEp}
then $J$ satisfies Hyp.~\ref{hyp:J_nonexp} (but not Hyp.~\ref{hyp:J_geo}) and it follows from Theorems~\ref{theo:speed_estimate_inf} and~\ref{theo:speed_estimate_sup} that the positions of the level sets 
$E_\lambda(t)$  are asymptotically algebraic and superlinear as $t\to+\infty,$ in the sense that for $\epsilon>0$, there is $\tilde\rho>f'(0)$ such that
\be
\lp\frac{f'(0)-\epsilon}{\beta}\rp^{1/\alpha}\ t^{1/\alpha}  \leq |x_\lambda(t)| \leq  \lp\frac{\tilde \rho}{\beta}\rp^{1/\alpha}\ t^{1/\alpha}\  \hbox{ for large } t,
\ee
where $x_\lambda(t)$ is any element of the level set $E_\lambda(t)$ (see Fig.~\ref{fig:xl_NEp}).

\begin{figure}[ht]
\begin{center}
 \subfigure[Accelerated propagation]{\includegraphics[width=17cm]{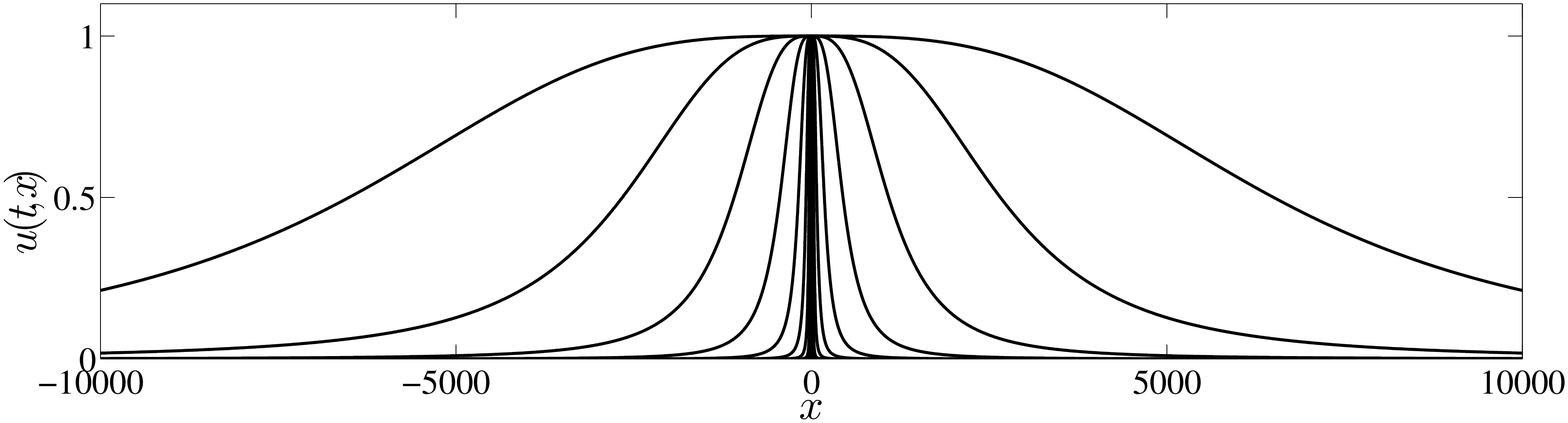}}
 \subfigure[Finite speed of propagation]{\includegraphics[width=17cm]{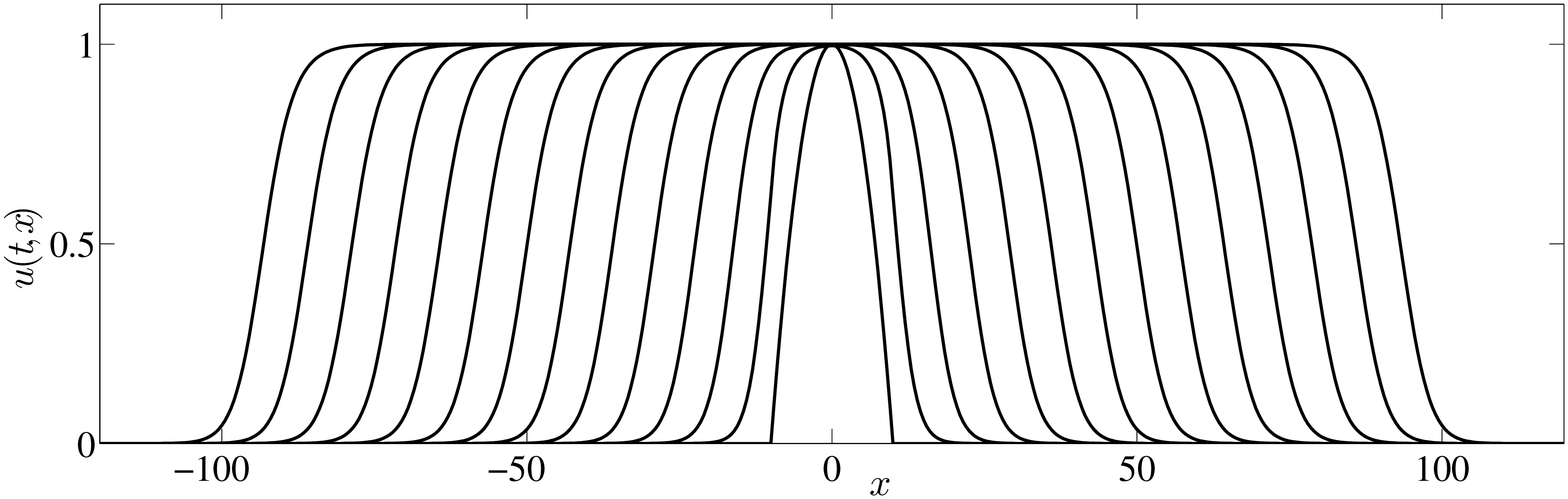}}
\end{center}
\caption{ The solution $u(t,x)$ of problem~\eqref{eq:1} at successive times $t =0,3,\dots,30,$ with $f(u)=u(1-u)$ and  $u_0(x)=\max((1-(x/10)^2),0):$ (a) with an exponentially 
unbounded kernel $J(x)=(1+|x|)^{-3};$
(b) with an exponentially bounded kernel
 $J(x)=(1/2)e^{-|x|}$.
}
\label{fig:u_geo_vs_exp}
\end{figure}

Next, if $J$ satisfies \eqref{hyp:J1} and decays algebraically for large $x$ as in~\eqref{exp:Geo},
then $J$ satisfies both Hyp.~\ref{hyp:J_nonexp} and~\ref{hyp:J_geo} and it follows from Theorems~\ref{theo:speed_estimate_inf} and~\ref{theo:speed_estimate_sup} that the position of the 
level sets
$E_\lambda(t)$ move exponentially fast as $t\to+\infty$ in the sense that, for any $\lambda\in(0,1)$ and $\epsilon>0$, there is $\tilde\rho>f'(0)$ such that
\be
 \ds \frac{f'(0)-\epsilon}{\alpha}\ t  \ \leq\  \ln{(|x_\lambda(t)|)}\ \leq \ \frac{\tilde \rho}{\alpha}\ t \ \hbox{ for large } t,
 \label{eq:estim_exple_geo}
\ee
for any $x_\lambda(t)$ in the level set $E_\lambda(t).$ The profile of the solution $u(t,x)$ of \eqref{eq:1} with an algebraically decreasing kernel is illustrated in Fig.~\ref{fig:u_geo_vs_exp} (a).

We mention that Cabr\'e and Roquejoffre~\cite{CabRoq09} just established comparable estimates
for the level sets of the solutions $u$ of equations of the type $u_t = Au + f(u),$ where $f$ is concave, $u_0$ is compactly supported or monotone one-sided compactly supported, 
and the operator $A$ is the generator of a Feller semi-group. A typical example is the fractional
Laplacian $A = -(-\Delta)^{\alpha}$ with $0 < \alpha < 1$: if $u$ is smooth enough and decays slowly to $0$ at infinty,
\[
\ds (-\Delta)^{\alpha}u(x)=c_\alpha \int_{\R}{\frac{u(y)-u(x)}{|x-y|^{1+2\alpha}}dy},\ \hbox{ for all } x\in\R,
\]
where $c_\alpha$ is choosen such that the symbol of $(-\Delta)^{\alpha}$ is $|\xi|^{2\alpha}.$
In this case, the asymptotic exponential spreading
of the level sets also follows from the algebraic decay of the kernel $J_\alpha(x)=|x|^{-(1+2\alpha)}$ associated with the operator $A.$ We can notice that it is not a particular case of
Theorems~\ref{theo:speed_estimate_inf} and~\ref{theo:speed_estimate_sup} since the kernel $J_\alpha$ is singular at $x=0.$

Lastly, let us consider the example of a function $J$ satisfying \eqref{hyp:J1},\eqref{hyp:J2} and such that $|J'/J|$ is not monotone as $|x|\rightarrow\infty$, e.g.
$$J(x)= C |x-\sin(x)|^{-\alpha}  \hbox{ for large } |x| 
\ \hbox{(with }\alpha >2).$$
Then $J$ does not satisfies Hyp. 1, but still satisfies Hyp. 2. Thus, we can apply Theorems~\ref{theo:speed_estimate_inf} 
and~\ref{theo:speed_estimate_sup} which lead to the same estimates as \eqref{eq:estim_exple_geo}.

In all above examples the positions of the level sets increase super-linearly with time. This illustrates the accelerating behavior of the solution of \eqref{eq:1} for exponentially unbounded kernels.
Coming back to Fig.~\ref{fig:u_geo_vs_exp} (a), we indeed observe that the distance between level sets of the same level tends to increase with time when time growths as $t_n=an$ 
with $a>0$ and when $J$ is exponentially unbounded, whereas it remains constant in the exponentially bounded case (Fig.~\ref{fig:u_geo_vs_exp} (b)). 
Moreover, in Fig.~\ref{fig:u_geo_vs_exp} (a), the solution tends to flatten as $t\to \infty$ \ie  the lower the level $\lambda,$ the faster the level sets $E_\lambda(t)$ propagate. 
In particular, this implies that the solution does not converge to a traveling wave solution. This is coherent with the fact that~\eqref{eq:1} does not admit traveling wave solutions 
when the kernel $J$ is exponentially unbounded~\cite{Yag09}.

\section{Proofs of the Theorems}
\subsection{Proof of Theorem~\ref{theo:infinite_speed}}
We begin with proving that for any $t\geq0,$ $\ds\liminf_{x\to\pm\infty}{u(t,x)}=0.$
Let us define $\bv(t,x)=v(t,x)e^{rt}$ for all $(t,x)\in[0,\infty)\times\R,$ where $\ds r=\sup_{s\in(0,1]}{(f(s)/s)}\geq f'(0)>0$ and $v$ satisfies the following problem:
\begin{equation}
  \left\{
           \begin{array}{l}
            v_t=J\ast v- v , \ t>0,\ x\in\R,\\
            v(0,x)=u_0(x),\ x\in\R.
           \end{array}
\right.
\label{eq:bv}
\end{equation}
Then $\bv$ verifies $\bv_t=J\ast \bv- \bv + r\bv$ on $(0,\infty)\times\R$ and $\bv(0,x)=u_0(x).$
From the maximum principle~\cite{Wei82,Yag09}, we get $0<u(t,x)\leq \bv(t,x)$ for all $(t,x)\in(0,\infty)\times\R.$

Moreover, since $u_0$ is compactly supported and the operator $u\mapsto J\ast u-u$ is Lipschitz-continuous on $L^\infty(\R)\cap L^1(\R),$ the Cauchy-Lipschitz theorem implies that the 
solution $t\mapsto v(t,\cdot)$ of problem~\eqref{eq:bv} belongs to $\C1([0,\infty),L^\infty(\R)\cap L^1(\R)).$ Integrating~\eqref{eq:bv} over $\R$ and using~\eqref{hyp:J1}, we get $\ds\|v(t)\|_{L^1(\R)}=\|u_0\|_{L^1(\R)}.$
This implies that $u(t)$ belongs to $L^1(\R)$ for all $t>0.$ Since $u(t,x)>0$ for any $t>0$ and $x\in\R,$ we have:
 \begin{equation}\label{eq:lim_u_infty}
\liminf_{x\to-\infty}{u(t, x)} =0 \hbox{ and } \liminf_{x\to+\infty}{u(t, x)} =0 \hbox{ for each }t\geq 0.
\end{equation}

Let us now prove that $u(t,0)\to1$ as $t\to\infty.$ Let $\tilde f:[0,1]\to\R$ satisfies~\eqref{hyp:f},~\eqref{hyp:f1} and such that $\tilde f\leq f$ in $[0,1],$ 
$\tilde f(s)\leq\tilde f'(0)s$ for all $s\in[0,1]$ and $\tilde f$ is a nonincreasing function in a neighborhood of $1.$
We denote by $\tilde u$ the solution of the Cauchy problem~\eqref{eq:1} with the nonlinearity $\tilde f.$ From the maximum principle $u\geq\tilde u$ on $[0,\infty)\times\R.$

Then, let us set $g_\epsilon(s)=\tilde f(s)-\epsilon s$ in $[0,1],$ where $\epsilon\in(0,1)$ is small enough such that $g_\epsilon'(0)>0.$ Set
$$\lambda_\epsilon=\sup\{s>0\ |\ g_\epsilon>0\hbox{ in } (0,s]\}<1.$$
One can choose $\epsilon>0$ small enough so that $g_\epsilon<0$ on $(\lambda_\epsilon,1].$ 
From~\eqref{hyp:J1} we know that there exists $A_\epsilon>0$ such that 
$$\ds D_\epsilon=\int_{-A_\epsilon}^{A_\epsilon}{J(y)dy}=1-\epsilon.$$

Let $\uv$ be the solution of the following Cauchy problem:
\begin{equation}\label{eq:J_supp_compact}\left\{
\begin{array}{l}
  \partial_t \uv=D_\epsilon(J_\epsilon\ast \uv- \uv) +g_\epsilon(\uv), \ t>0, \ x\in\R, \\
  \uv(0,x)=u_0(x), \ x\in\R,
\end{array}
\right.
\end{equation}
where $J_\epsilon$ is a compactly supported kernel defined by:
\[
J_\epsilon(x)=\frac{J(x)}{\ds \int_{-A_\epsilon}^{A_\epsilon}{J(y)dy}}\mathds{1}_{[-A_\epsilon,A_\epsilon]}(x).
\]
We have 
$$\partial_t \uv=(J\times \mathds{1}_{[-A_\epsilon,A_\epsilon]})\ast \uv- \uv +\tilde f(\uv)\le J \ast \uv- \uv +\tilde f(\uv).$$
The maximum principle implies that 
$0\leq\uv(t, x)\leq \tilde u(t,x) \leq u(t, x)$ ($\leq 1$) for all $t \geq 0$ and $x \in\R$.
From Theorem~$3.2$ in~\cite{LutPasLew05}, we know that $\uv$ propagates with a finite speed $c^*_\epsilon>0$ \ie for all $c\in(0,c^*_\epsilon),$

\begin{equation}\label{eq:prop_speed}
 \sup_{|x|\leq ct}{|\uv(t,x)-\lambda_\epsilon|}\to0 \hbox{ as } t\to+\infty.
\end{equation}
In particular, we have:
\[
 \lim_{t\to\infty}{\uv(t,0)}=\lambda_\epsilon \leq \lim_{t\to\infty}{u(t,0)}\leq 1.
\]
Since $\lambda_\epsilon\to1$ as $\epsilon\to0,$ we get
\[
 u(t,0)\to1 \hbox{ as }t\to\infty.
\]
It then follows that for any $\lambda\in(0, 1)$ there exists a time $t_\lambda \geq 0$ such that
\begin{equation}\label{eq:lim_u_0}
u(t,0)>\lambda \hbox{ for all }t \geq t_\lambda.
\end{equation}
Since the functions $x \mapsto u(t, x)$ are continuous for all $t>0,$ one concludes from~\eqref{eq:lim_u_infty} and~\eqref{eq:lim_u_0} that $E_\lambda(t)$ is a non-empty set for all $t \geq t_\lambda.$

Let us now prove~\eqref{eq:E}. From~\cite{LutPasLew05}, we know that the propagation speed $c^*_\epsilon$ is the minimal speed of traveling wave solutions of problem~\eqref{eq:J_supp_compact}. 
This speed verifies~\cite{CarChm04,Sch80}:
\be
\ds c^{*}_\epsilon=\min_{\eta>0}{\left\{ \ds\frac{1}{\eta}\lp D_\epsilon\int_{\R}{J_\epsilon(z)e^{\eta z}dz} -1 +\tilde f'(0)\rp\right\}}.
\label{eq:c}
\ee
Let $A>0,$ $\lambda\in(0,1),$ $t\geq t_\lambda$ and  $x_\lambda(t)\in E_\lambda(t).$ Since $J$ is exponentially unbounded in the sense of~\eqref{hyp:J3}, 
$c^*_\epsilon\to \infty$ as $\epsilon\to0.$ Let us choose $\epsilon>0$ small enough such that $\lambda_\epsilon>\lambda$ and $c^*_\epsilon>A.$ Then, it follows from~\eqref{eq:prop_speed} that $|x_\lambda(t)|\geq At$
for $t$ large enough. Since this is true for any $A>0$ and any $x_\lambda(t)\in E_\lambda(t)$ we get~\eqref{eq:E_lin} and~\eqref{eq:E}. \carre

\subsection{Proof of Theorem~\ref{theo:speed_estimate_inf}}

This section is devoted to the proof of a lower bound for $\min\{E_\lambda(t)\cap (0,+\infty)\}$ (resp. $-\max\{E_\lambda(t)\cap (-\infty,0)\}$).
The proof is divided into three parts. We begin with showing that the solution of~\eqref{eq:1} at time $t=1$ is larger than some multiple of $J.$ Then, we construct an 
appropriate subsolution of~\eqref{eq:1} which enables us to prove the lower bound for small values of $\lambda$. Lastly, we show that the lower bound remains true for 
any value of $\lambda\in (0,1).$

More precisely, let us fix $\lambda \in (0, 1)$  and $\epsilon\in(0,f'(0)).$ We claim that
\begin{equation}
 E_\lambda(t)\subset J^{-1}\left\{\left(0,e^{-(f'(0)-\epsilon)t}\right]\right\}
\label{eq:ineq_inf_E}
\end{equation}
for $t$ large enough.
\

\paragraph*{\underline{Step 1:} $u(1,\cdot)$ is bounded from below by a multiple of $J(\cdot)$}

Let us define 
$$v(t,x)=(u_0(x)+t(J\ast u_0)(x))e^{-t}.$$
Then it is easy to see that $v$ is a subsolution of the following linear Cauchy problem:
\be
\left\{
\begin{array}{l}
 \uu_t=J\ast \uu-\uu, \ t>0, \ x\in \R, \\
 \uu(0,x)=u_0(x), \ x\in \R.
\end{array}
\right.
 \label{eq:lin_inf}
\ee

Indeed, $v(0,x)=u_0(x)$ and  for all $(t,x)\in[0,\infty)\times\R$:
\[
 \ds v_t-J\ast v+v=-t e^{-t} \lp J\ast \lp J\ast u_0\rp\rp(x)\leq 0.
\]
Since $\uu$ is also a subsolution of the equation~\eqref{eq:1} verified by $u,$ we get:
\[
 u(t,x )\geq v(t,x) \hbox{ for all } (t,x)\in[0,\infty)\times\R.
\]

Moreover,  $v(1,x)=(u_0(x)+(J\ast u_0)(x))e^{ -1}$ for all $x\in\R$ which implies that is exists $C\in(0,1)$ such that $v(1,x)\geq CJ(x)$ for all $x\in\R.$
Finally,
\begin{equation}\label{eq:ineq_uvJ}
  u(1,x )\geq v(1,x)\geq CJ(x) \hbox{ for all } x\in\R.
\end{equation}

\

\paragraph*{\underline{Step 2:} Proof of~\eqref{eq:ineq_inf_E} for small values of $\lambda$}

We recall that there exist $\delta > 0,$ $s_0 \in (0, 1)$ and $M \geq 0$ such that $f (s) \geq f'(0) s - M s^{1+\delta}$ for all $s \in [0, s_0 ].$

Define $\rho_1>0$ by
\begin{equation}\label{eq:rho1}
 \rho_1=f'(0)-\epsilon/2.
\end{equation}

Since $J$ satisfies~\eqref{hyp:J1} and~\eqref{hyp:J2}, we can choose $\xi_1>0$ such that for all $|x|\geq \xi_1,$
\begin{equation}
 \left|J'(x)\right|\leq (\epsilon'/2) \times J(x),
\end{equation}
where $\epsilon'>0$ satisfies
\[
 \epsilon'\int_{0}^\infty{J(z)zdz}\leq\epsilon/2.
\]
Let us set:
\begin{equation}\label{eq:s1_Ba}
\ds\kappa=\inf_{(-\xi_1,\xi_1)}{CJ}=CJ(\xi_1)>0, \ s_1=\min{(s_0,\kappa)},
\end{equation}
and
\begin{equation}\label{eq:s1_Bb}
 B=\max{\left(s_1^{-\delta},\ds \frac{M}{\rho_1\delta},
\ds\frac{\ds M\left(\ds\frac{\delta}{1+\delta}\right)^{\delta}}{(1+\delta)(f'(0)-\epsilon/2)} \right)}>0.
\end{equation}
Let $g$ be the function defined in $[0,\infty)$ by
\[
 g(s)=s-Bs^{1+\delta}.
\]
We observe that
\[
 g(s)\leq0 \hbox{ for all } s\geq s_1 \hbox{ and } g(s)\leq s_1 \hbox{ for all } s\geq0.
\]
Moreover, let $0<s_2<s_1$ be such that $g'(s_2)=0$ and
\begin{equation}\label{eq:lambda_2}
\ds\lambda_2=g(s_2)=\max_{s\in[0,s_0]}{g(s)}=\frac{\delta}{1+\delta}\left((1+\delta)B\right)^{-1/\delta}.
\end{equation}
Let $\xi_0(t)>0$ be such that:
\be
\label{eq:defxi0}
 C\, J(\xi_0(t))e^{\rho_1 t}=s_2 \hbox{ for all }t\ge0.
\ee
We can notice that for all $t\geq0$, $\xi_0(t)\geq\xi_1$ and that $\xi_0(t)$ is continuous and increasing in $t\geq0,$ since $J$ is continuous and decreasing in $[0,\infty).$

Then, let us define $\uu$ as follows:
\begin{equation}\label{eq:def_uu}
 \ds \uu(t,x)=\left\{ \begin{array}{ll}
                 \ds g(CJ(x)e^{\rho_1 t})                       & \hbox{for all } |x|> \xi_0(t) \\
                 \ds \lambda_2=g(s_2)=g(CJ(\xi_0(t))e^{\rho_1 t}) & \hbox{for all } |x|\leq\xi_0(t)
                 \end{array} \right.
\hbox{ for all } t\geq0.
\end{equation}
Observe that $0<CJ(x)e^{\rho_1 t}\leq CJ(\xi_0(t))e^{\rho_1t}=s_2< s_1$ when $|x|\geq\xi_0(t),$ whence $\uu(t,x)>0$ for all $t\geq0$ and $x\in\R.$ 
Let us check that $\uu$ is a sub-solution of~\eqref{eq:1}. Since $J(x)$ is nonincreasing with respect to $|x|$ and $u$ satisfies~\eqref{eq:ineq_uvJ}, we have
\be
\label{eq:ineg_uu0}
\begin{array}{l}
 \ds \uu(0,x)=CJ(x)-B(CJ(x))^{1+\delta}\leq CJ(x)\leq u(1,x) \hbox{ for all }|x|>\xi_0(0) \\
 \ds \uu(0,x)=\lambda_2=g(CJ(\xi_0(0)))\leq CJ(x)\leq u(1,x) \hbox{ for all }|x|\leq\xi_0(0).
\end{array}
\ee 

Then, let us check that $\uu$ is a subsolution of the equation satisfied by $u$ in the region where $\uu<\lambda_2$. Let $(t,x)$ be any point
in $[0,\infty)\times\R$ such that $\uu(t,x)<\lambda_2.$As already emphasized, one has $0<CJ(x)e^{\rho_1 t}<s_1,$ whence $CJ(x)<s_1$ and
$|x|\geq\xi_1$ from~\eqref{eq:s1_Ba}. Furthermore,
\begin{equation}
 0<\uu(t,x)<CJ(x)e^{\rho_1 t}<s_1\leq s_0<1.
\end{equation}
Thus, since $f$ satisfies~\eqref{hyp:f1}, we get
\begin{equation}\label{ineq:f}
 f(\uu(t,x))\geq f'(0)\left( CJ(x)e^{\rho_1 t} -B (CJ(x))^{1+\delta}e^{\rho_1(1+\delta) t}\right) -M(CJ(x))^{1+\delta}e^{\rho_1(1+\delta) t}.
\end{equation}

Let us now show that $J\ast\uu-\uu\geq -(\epsilon/2)\uu,$ for all $t\in[0,\infty)$ and $|x|>\xi_0(t).$
Let $t\in[0,\infty)$ and $$x>\xi_0(t)>0.$$ 
Remember that $J$ is decreasing on $[0,+\infty)$ and tha $g$ is decreasing in $[0,s_2]$ (which implies that $\uu(t,y)$ is nonincreasing with respect to $|y|$). Then 
\begin{equation}\label{ineq:Ju_u}
\begin{array}{rcl}
 \ds (J\ast\uu)(t,x)-\uu(t,x)& = & \ds\int_{\R}{J(x-y)(\uu(t,y)-\uu(t,x))dy}\\ [0.2cm]
                  & \geq & \ds \int_{|y|>x}{J(x-y)(\uu(t,y)-\uu(t,x))dy} \\ [0.2cm]
                  & \geq & \ds \int_{-\infty}^{-x}{J(x-y)(\uu(t,y)-\uu(t,x))dy} +\int_{x}^\infty{J(x-y)(\uu(t,y)-\uu(t,x))dy} \\ [0.2cm]
                  & \geq & \ds \int_{x}^\infty{(J(x-y)+J(x+y))(\uu(t,y)-\uu(t,x))dy}.
\end{array}
\end{equation}
Observe that for all $y>x,$
\[
 0\geq\uu(t,y)-\uu(t,x)=\int_{x}^{y}{\partial_x{\uu}(t,s)ds}.
\]
Furthermore, for all $t>0$ and $s>\xi_0(t)$
\[
 \partial_x{\uu}(t,s)=CJ'(s)e^{\rho_1 t}g'\left(CJ(s)e^{\rho_1 t}\right)= CJ'(s)e^{\rho_1 t}\left( 1-(1+\delta)B (CJ(s))^{\delta}e^{\rho_1\delta t}\right).
\]
Since $s>\xi_0(t)\geq \xi_1:$
\[
\begin{array}{rcl}
 |\partial_x{\uu}(t,s)| & \leq & (\epsilon'/2) CJ(s)e^{\rho_1 t}g'\left(CJ(s)e^{\rho_1 t}\right) \\
                        & \leq & (\epsilon'/2) CJ(s)e^{\rho_1 t}\left( 1-B (CJ(s))^{\delta}e^{\rho_1\delta t}\right)= (\epsilon'/2) \uu(t,s) \\
                        & \leq & (\epsilon'/2)\uu(t,x).
\end{array}
\]
Finally,  for all $y\geq x>\xi_0(t)$
\begin{equation}\label{ineq:uu_epsilon}
 \uu(t,y)-\uu(t,x)\geq -(\epsilon'/2)(y-x)\uu(t,x).
\end{equation}
Then, from equation~\eqref{ineq:Ju_u} and~\eqref{ineq:uu_epsilon}, we get
\begin{equation}\label{ineq:Juu_uu}
\begin{array}{rcl}
 \ds (J\ast\uu)(t,x)-\uu(t,x)& \geq & \ds \int_{x}^\infty{(J(x-y)+J(x+y))(\uu(t,y)-\uu(t,x))dy} \\ [0.2cm]
                           & \geq & \ds -(\epsilon'/2)\left(\int_{x}^\infty{J(y-x)(y-x)dy} + \int_{x}^\infty{J(x+y)(y-x)dy}\right)\uu(t,x) \\ [0.2cm]
                           & \geq & \ds -(\epsilon'/2)\left(\int_{0}^\infty{J(z)zdz} + \int_{0}^\infty{J(2x+z)zdz}\right)\uu(t,x) \\ [0.2cm]
                           & \geq & \ds -\epsilon'\lp\int_{0}^\infty{J(z)zdz}\rp\ \uu(t,x) \\ [0.2cm]
                           & \geq & \ds -(\epsilon/2)\, \uu(t,x).
\end{array}
\end{equation}
The same property holds for $x<-\xi_0(t)$ by symmetry of $J$ and $\uu$ with respect to $x.$
It follows from~\eqref{eq:rho1},~\eqref{eq:s1_Bb},~\eqref{ineq:f} and~\eqref{ineq:Juu_uu} that, for all $t\geq0$ and $|x|\geq\xi_0(t)$
\begin{equation}\label{eq:uu_sup}
\begin{array}{l}
\ds \uu_t(t,x) -\left(J\ast \uu(t,x)- \uu (t,x)\right) - f(\uu(t,x))\\ [0.2cm]
 \begin{array}{cl}
 \leq & \ds \rho_1 CJ(x)e^{\rho_1 t} -\rho_1(1+\delta)B(CJ(x))^{1+\delta}e^{\rho_1(1+\delta)t} + \epsilon/2 \, \uu(t,x) \\ [0.2cm]
      & \ds -f'(0)\left( CJ(x)e^{\rho_1 t} -B (CJ(x))^{1+\delta}e^{\rho_1(1+\delta) t}\right) + M(CJ(x))^{1+\delta}e^{\rho_1(1+\delta) t}  \\[0.2cm]

 \leq & \ds \left(\rho_1 -f'(0) + \epsilon/2 \right)CJ(x)e^{\rho_1 t} +\left(M -B[\rho_1(1+\delta)-f'(0)+\epsilon/2]\right)(CJ(x))^{1+\delta}e^{\rho_1(1+\delta)t} \\ [0.2cm]
\leq & (M -B \, \delta \, \rho_1)(CJ(x))^{1+\delta}e^{\rho_1(1+\delta)t}\\
 \leq & 0.
 \end{array}
\end{array}
\end{equation}

Let us now check that $\uu$ is a subsolution of the equation satisfied by $u$, in the region where $\uu=\lambda_2$. Let $(t,x)$ be any point
in $[0,\infty)\times\R$ such that $\uu(t,x)=\lambda_2.$ The same arguments as above imply that for any point $(t,x)$ satisfying $\uu(t,x)=\lambda_2,$ \ie for all
$(t,x)\in[0,\infty)\times[-\xi_0(t),\xi_0(t)],$ we get
\begin{equation}\label{ineq:Ju_lambda}
\begin{array}{rcl}
 \ds (J\ast\uu)(t,x)-\uu(t,x)& = & \ds \int_{\R}{J(x-y)(\uu(t,y)-\lambda_2)dy} \\ [0.2cm]
                           & = & \ds \int_{|y|\geq\xi_0(t)}{J(x-y)(\uu(t,y)-\lambda_2)dy} \\ [0.2cm]
                           & = & \ds \int_{\xi_0(t)}^\infty{(J(x-y)+J(x+y))(\uu(t,y)-\uu(t,\xi_0(t)))dy} \\ [0.2cm]
                           & \geq & \ds -(\epsilon'/2)\left(\int_{\xi_0(t)}^\infty{J(y-x)(y-\xi_0(t))dy} + \int_{\xi_0(t)}^\infty{J(x+y)(y-\xi_0(t))dy}\right)\lambda_2 \\ [0.2cm]
                           & \geq & \ds -(\epsilon'/2)\left(\int_{0}^\infty{J(z+\xi_0(t)-x)zdz} + \int_{0}^\infty{J(z+x+\xi_0(t))zdz}\right)\lambda_2 \\ [0.2cm]
                           & \geq & \ds -\epsilon'\lp\int_{0}^\infty{J(z)zdz}\rp\ \lambda_2 \\ [0.2cm]
                           & \geq & \ds -(\epsilon/2)\, \lambda_2.
\end{array}
\end{equation}
For all $x\in [0,\xi_0(t))$ it follows from~\eqref{eq:def_uu} that $\uu_t(t,x)=0.$ Let us show that this is also true when $x=\xi_0(t).$  As already noticed, $\xi_0(t)$ is  
an increasing function of $t.$ Thus, for all $h>0$ $\xi_0(t)<\xi_0(t+h),$ which implies that $\uu(t+h,\xi_0(t))=\lambda_2=\uu(t,\xi_0(t)).$ As a consequence,
\[
 \ds\lim_{h\to0,h>0}{\ds\frac{\uu(t+h,\xi_0(t))-\uu(t,\xi_0(t))}{h}}=0.
\]
Moreover,
\begin{equation*}
\ds\lim_{h\to0,h>0}{\ds\frac{\uu(t,\xi_0(t))-\uu(t-h,\xi_0(t))}{h}} =  \ds\lim_{h\to0,h>0}{\ds\frac{\lambda_2-g(CJ(\xi_0(t))e^{\rho_1 (t-h)})}{h}}.
\end{equation*}
For all $h>0$ small enough and using the definition~\eqref{eq:defxi0} of $\xi_0(t)$, we obtain:
\[
g(CJ(\xi_0(t))e^{\rho_1 (t-h)})=\lambda_2-\rho_1 CJ(\xi_0(t))g'(s_2)h+o(h)=\lambda_2+o(h).
\]
This implies that
\begin{equation}\label{ineq:uu_t}
 \uu_t(t,\xi_0(t))  =  \ds\lim_{h\to0}{\frac{\uu(t,\xi_0(t))-\uu(t-h,\xi_0(t))}{h}}=0.
\end{equation}
Since $\uu_t(t,x)=0$ for all $x\in [0,\xi_0(t))$, the above equality and the symmetry of the problem imply that:
\be
\uu_t(t,x)=0 \hbox{ for all }|x|\le \xi_0(t).
\label{eq:derivu}
\ee
It follows from~\eqref{eq:s1_Bb},~\eqref{ineq:Ju_lambda} and~\eqref{eq:derivu} that for all $t\geq0$ and $|x|\leq\xi_0(t),$
\begin{equation}
\begin{array}{l}
\ds \uu_t(t,x) -\left(J\ast \uu(t,x)- \uu (t,x)\right) - f(\uu(t,x))\\ [0.2cm]
 \begin{array}{cl}
 \leq & \ds (\epsilon/2) \, \lambda_2 -\left(f'(0)\lambda_2 -M \lambda_2^{1+\delta}\right) \\[0.2cm]

 \leq & \ds -\lambda_2(f'(0)-\epsilon/2)\left( 1 -\frac{M}{(1+\delta)B(f'(0)-\epsilon/2)}\left(\frac{\delta}{1+\delta}\right)^{\delta} \right) \\ [0.2cm]

 \leq & 0.
 \end{array}
\end{array}
\end{equation}

Using the above inequality together with~\eqref{eq:ineg_uu0} and ~\eqref{eq:uu_sup}, the maximum principle implies that
\begin{equation}\label{ineq:uu}
 \uu(t-1,x)\leq u(t,x) \hbox{ for all } t\geq1 \hbox{ and } x\in\R.
\end{equation}

Fix now any real number $\omega$ small enough so that:
\[
 0<\omega<s_2 
\]
This real number $\omega$ does not depends on $\lambda$  but depends on $\epsilon,$ as well as on $J$ and $f.$ Remember that $t_\omega\geq0$ is
such that $E_\omega(t)$ is a non-empty set for all $t\geq t_\omega.$ Since $J$ is continuous and decreasing on $[0,+\infty),$ there exists then a time 
$\overline{t}_\omega\geq\max{(t_\omega,1)}$
such that for all $t\geq \overline{t}_\omega,$ it exists $y_\omega(t)\in(\xi_1,\infty)$ such that 
$$CJ(y_\omega(t))e^{\rho_1(t-1)}=\omega.$$
Furthermore, the function $\ds y_\omega(t):\left[\overline{t}_\omega,\infty\right)\to[\xi_1,\infty)$ is increasing and continuous.

Lastly, let $\Omega$ be the open set defined by
\[
 \Omega=\{ (t,x)\in(\overline{t}_\omega,+ \infty)\times\R, \ |x|<y_\omega(t) \}.
\]
We claim that $\ds\inf_\Omega{ u }> 0.$ Indeed, if $(t,x)\in\Omega$ is such that $CJ(x)e^{\rho_1(t-1)}\geq s_2,$ then $|x|\leq\xi_0(t-1)$ and 
$$\uu(t-1,x)=\lambda_2=g(s_2)>g(\omega)>0.$$ 
Otherwise, $(t,x)$ is such that
$\omega<CJ(x)e^{\rho_1(t-1)}<s_2,$ whence $|x|>\xi_0(t-1)$ and 
$$\uu(t-1,x)=g(CJ(x)e^{\rho_1(t_1)})\geq g(\omega)>0.$$ 
Finally equation~\eqref{ineq:uu} implies that
\begin{equation}\label{eq:ineq_ut}
 u(t, x) \geq g(\omega)>0 \hbox{ for all } (t, x) \in\Omega.
\end{equation}
Thus, setting $\theta=g(\omega),$ we get that if $\lambda\in (0, \theta)$ and if $x\in E_\lambda (t)$ for $t \geq\max(t_\lambda , \tom ),$ then
\[
 |x|\geq y_\omega(t)\geq\xi_1\geq\xi_0(0)
\]
Since  $\rho_1 = f'(0) - \epsilon/2 >f'(0)-\epsilon,$ there exists then a time $T_{\lambda,\epsilon} \geq\max(t_\lambda , \tom )$ such
that,
\begin{equation}\label{ineq:inf_E-}
 \forall t>T_{\lambda,\epsilon},\ \forall x\in E_\lambda(t), \  J(x)\leq J(y_\omega(t))=\frac{\omega e^{\rho_1} }{C}e^{-\rho_1 t}\leq  e^{-(f'(0)-\epsilon)t}.
\end{equation}
This proves~\eqref{eq:ineq_inf_E} for $\lambda\in(0, \theta)$.

\

\paragraph*{\underline{Step 3:} Proof of \eqref{eq:ineq_inf_E} for any $\lambda\in(0,1)$}

Assume that $\lambda\in(0,1).$ Let $\uuto$ be the function defined by
\[
 \uuto(x)=\left\{
                 \begin{array}{ll}
                 \theta(1 -|x|)    & \hbox{if } |x|\leq1, \\
                 0         & \hbox{if }  |x|>1,
                \end{array}
           \right.
\]
where $\theta=g(\omega)$ is given in Step~$2.$ Let us set  $\tilde f(s)=f(s)-(1-D_\lambda)s$ for all $s\in[0,1]$ with 
$$\ds D_\lambda=\int_{-\xi_1}^{\xi_1}{J(y)dy}\in(0,1),$$
where $\xi_1$ is chosen large enough such that $\tilde f'(0)>0$ and $\tilde f(s)>0$ in $(0,\lambda].$

We consider the solution $\uut$  of the Cauchy problem:
\begin{equation}\left\{
\begin{array}{l}
  \partial_t\uut=D_\lambda\left(J_\lambda\ast \uut- \uut\right) +\tilde f(\uut), \ t>0, \ x\in\R, \\
  \uut(0,x)=\uuto(x), \ x\in\R,
\end{array}
\right.
\end{equation}
where $J_\lambda$ is a compactly supported kernel defined by:
\[
 \ds J_\lambda(x)=\frac{J(x)}{\ds \int_{-\xi_1}^{\xi_1}{J(y)dy}}\mathds{1}_{\ds[-\xi_1,\xi_1]}(x), \ x\in\R.
\]
It follows from~\eqref{eq:ineq_ut} that
\[
 \forall T\geq\tom, \ \forall |\xi|\leq y_\omega(T)-1, \ \forall \ x\in \R,\ u(T,x)\geq \uuto(x- \xi)),
\]
whence
\begin{equation}\label{ineq:uT}
 \forall T\geq\tom, \ \forall |\xi|\leq y_\omega(T)-1, \  \forall t\geq0,\ \forall x\in\R, \ u(T+t,x)\geq \uut(t,x-\xi)
\end{equation}
from the maximum principle. Indeed, we have for all $(t,x)\in(0,\infty)\times\R$:
\[
\begin{array}{l}
\ds \partial_t\uut(t,x) - \left(J\ast \uut(t,x) - \uut(t,x)\right) - f(\uut(t,x)) \\ [0.2cm]
 \begin{array}{cl}
 \leq & \ds \partial_t\uut(t,x) - D_\lambda J_\lambda\ast \uut(t,x) + \uut(t,x) - f(\uut(t,x)) \\[0.2cm]

 \leq & \ds \partial_t\uut(t,x) - D_\lambda\left(J_\lambda\ast \uut(t,x)- \uut(t,x)\right) - \tilde f(\uut(t,x)) \\ [0.2cm]

 \leq & 0.
 \end{array}
\end{array}
\]
Moreover, we know from Theorem~$3.2$ in~\cite{LutPasLew05} that there exists $c^*_\lambda>0$ such that
\[
 \liminf_{t\to\infty}{\inf_{|x|<c^*_\lambda t}{\uut(t,x)}}=\sup{\{s>0\ |\ \tilde f>0 \hbox{ in } (0,s) \}}>\lambda.
\]
In particular, there exists $T_\lambda\geq0$  such that $\uut(T_\lambda,x)>\lambda$. Therefore,
\eqref{ineq:uT} implies that
\begin{equation*}
 \forall T\geq\tom, \ \forall |x|\leq y_\omega(T)-1, \ u(T+T_\lambda,x)>\lambda.
\end{equation*}
As a consequence, there exists $\underline{T}_{\lambda,\epsilon}\geq\max{(\tom + T_\lambda , t_\lambda )}$ such that for all $t \geq \underline{T}_{\lambda,\epsilon}$ and
for all $x \in E_\lambda (t),$ one has $|x| > y_\omega(t - T_\lambda )-1$ and
\begin{equation}\label{ineq:inf_E_+}
 J(x)\leq J\left( y_\omega\left(t -T_\lambda \right)-1\right)
= \frac{\omega e^{\rho_1}}{C}e^{-\rho_1(t - T_\lambda )}\times\frac{J\left( y_\omega\left(t -T_\lambda \right)-1\right)}{J\left( y_\omega\left(t -T_\lambda \right)\right)}
\leq  e^{-(f'(0)-\epsilon)t},
\end{equation}
since $y_\omega(t-T_\lambda)\to\infty$ as $t\to\infty$ and $\ds\frac{J(s-1)}{J(s)}\to1$ as $s\to\infty,$ from~\eqref{hyp:J2}.
This implies~\eqref{eq:ineq_inf_E} and completes the proof of Theorem~\ref{theo:speed_estimate_inf}. \carre

\subsection{Proof of Theorem~\ref{theo:speed_estimate_sup}}

In this section, we prove an upper bound for $\max\{E_\lambda(t)\cap (0,+\infty)\}$ (resp. $-\min\{E_\lambda(t)\cap (-\infty,0)\}$). The proof of this upper bound is based on the construction of suitable
supersolutions of~\eqref{eq:1}. The construction of such supersolutions strongly relies on Hypotheses~\ref{hyp:J_nonexp} and~\ref{hyp:J_geo}.

We shall prove that there exists $\rho>0$ such that, for any $\lambda\in(0,1),$
\begin{equation}
 E_\lambda(t)\subset J^{-1}\left\{[ e^{-\rho t} ,J(0)]\right\} \hbox{ for large }t.
\label{eq:ineq_sup_E}
\end{equation}

Since $f\in \C1([0,1]),$ the ``per capita growth rate" $f(s)/s$ is bounded from above by $\ds r=\sup_{s\in(0,1]}{\lp f(s)/s\rp}>0.$

\

\paragraph*{ \underline{Proof of \eqref{eq:ineq_sup_E} under Hypothesis \ref{hyp:J_nonexp}} }\label{sec:hyp_J_nonexp}

Assume that $J$ satisfies Hyp.~\ref{hyp:J_nonexp}. Then $J'/J$ is negative and nondecreasing on $[\sigma,+\infty)$. Therefore $\ln(J)$ is a nonincreasing convex function on 
$[\sigma,+\infty).$ It is then possible to define a function $\varphi: \, [0,\infty)\to [0,\infty)$ and $\tau\in[0,\sigma]$ such that $\varphi$ is nondecreasing and concave and 
\be
\varphi(0)=0, \, \varphi(+\infty)=+\infty, \hbox{ and }\varphi(x)=(\epsilon_0-1)\ln(J(x+\tau)) \hbox{ on }[\sigma-\tau,+\infty).
\label{eq:prop_phi}
\ee

By concavity, we have the following property, for all $y\geq x\geq0$:
\be
 \varphi(y)-\varphi(x)\leq \varphi(y-x).
\label{eq:ineg_phi}
\ee
Thus, we claim that:
\begin{equation}
 \label{eq:ineg_phi2}
\forall \, x\geq0, \, \forall \, y\in\R, \  \varphi(x)-\varphi(|x-y|)\leq \varphi(|y|).
\end{equation}
Even if it means increasing $\sigma,$ one can assume without loss of generality, that $J<1$ on $[\sigma,\infty).$
Indeed, if $y\leq0$ then $\varphi(x)-\varphi(|x-y|)\leq 0 \le \varphi(|y|).$ If $y>0,$ since $\varphi$ is nondecreasing and from~\eqref{eq:ineg_phi}, we have 
$\varphi(|y-x|)\geq\varphi(\max(x,y))-\varphi(\min(x,y))\geq \varphi(x)-\varphi(y).$ Notice that~\eqref{eq:ineg_phi2} implies immediately that 
\begin{equation}
\label{eq:ineg_phi3} 
\forall x\in\R,\  \forall y\in\R,\ \varphi(|x|)-\varphi(|x-y|)\leq \varphi(|y|).
\end{equation}

Let us define \[
 \phi(x)= e^{-\varphi(|x|)} \hbox{ for all }x\in \R.
\]
Using \eqref{eq:ineg_phi3}, we get:
\begin{equation}
\label{eq:J*phi/phi}
\begin{array}{rcl}
 \ds\frac{(J\ast \phi)(x)}{\phi(x)} & = & \ds\int_{\R}{\frac{\phi(x-y)}{\phi(x)}J(y)dy} \\ [0.4cm]
                                  & = & \ds \int_{\R}{J(y)e^{\varphi(|x|)-\varphi(|x-y|)}dy}   \\ [0.4cm]
                                  & \leq & \ds \int_{\R}{J(y)e^{\varphi(|y|)}dy}=\int_{\R}{J(y)/\phi(y)}dy.
\end{array}
\end{equation}
Moreover,
\begin{equation}
\begin{array}{rcl}
 \ds \int_{\R}{J(y)/\phi(y)}dy = \int_{\R}{J(y)e^{\varphi(|y|)}dy} 
                         & = &  \ds \int_{|y|<\sigma}{J(y)e^{\varphi(|y|)}dy}+\int_{|y|\ge \sigma }{J(y)e^{\varphi(|y|)}dy}\\ [0.4cm]
                         & = &  \ds \int_{|y|<\sigma}{J(y)e^{\varphi(|y|)}dy}+\int_{|y|\ge \sigma }{\lp\frac{J(y+\tau)}{J(y)}\rp^{\epsilon_0-1}J(y)^{\epsilon_0}}dy<\infty,  
\end{array}
\end{equation}
from Hyp.~\ref{hyp:J_nonexp} and since $\ds\frac{J(y+\tau)}{J(y)}\to1$ as $y\to\pm\infty$ from~\eqref{hyp:J2}.

Since $u_0$ is compactly supported, there exists $\sigma_1>0$ such that $\phi(x)\geq u_0(x),$ for all $|x|\geq\sigma_1.$
Finally, set 
$$\ds \rho_0=\max{ \left\{\int_{\R}{J(y)e^{\varphi(|y|)}dy},1\right\}}-1+r,$$ 
and define $\bu$ as follows:
\[
\forall (t,x)\in[0,+\infty)\times\R, \,\ds\overline{u}(t,x)=\min{\ds\left(\frac{\phi(x)}{\phi(\sigma_1)}e^{\rho_0 t},1\right)}.
\]

Observe that $u_0(x)\leq\bu(0,x)$ for all $x\in\R.$ Let us now check that $\bu$ is a supersolution of the equation~\eqref{eq:1} satisfied by $u.$
Since $u \leq 1,$
it is enough to check that $\bu$ is a supersolution of~\eqref{eq:1} whenever $\bu < 1.$ Note that since $\phi(x)$ is nonincreasing with respect to $|x|$, $\bu(t,x)< 1$ implies that $|x|>\sigma_1.$ 
Assume that $(t, x) \in[0, +\infty) \times[\sigma_1,\infty)$ and $\bu(t, x) < 1,$ then it follows from~\eqref{eq:J*phi/phi} that
\begin{equation*}
\begin{array}{rcl}
 \ds (J\ast\bu)(t,x) & \leq & \ds  (J\ast \phi)(x)\frac{e^{\rho_0 t}}{\phi(\sigma_1)} \\ [0.4cm]
                   & \leq & \ds  \frac{\phi(x)}{\phi(\sigma_1)}e^{\rho_0 t} \, \int_{\R}{J(y)/\phi(y)}dy \\ [0.4cm]
                   & = & \ds \bu(t,x)\, \int_{\R}{J(y)/\phi(y)}dy.
\end{array}.
\end{equation*}
This implies that for all $(t,x)$ such that $\bu(t,x)<1$
\begin{equation}
\begin{array}{l}
\bu_t(t,x)-(J\ast \bu)(t,x)+\bu(t,x)-f(\bu(t,x))\\ [0.2cm]
 \begin{array}{rl}
 \geq & \ds \rho_0 \bu(t,x) - (J\ast\bu)(t,x) +\bu(t,x) - r \bu(t,x) \\ [0.2cm]

 \geq & \ds (\rho_0 -  \int_{\R}{J(y)/\phi(y)}dy +1 -r )\bu(t,x)\\ [0.2cm]

 \geq & 0,
 \end{array}

\end{array}
\end{equation}
from the definition of $\rho_0.$
The parabolic maximum principle~\cite{Wei82,Yag09} then implies that:
\[
u(t,x)\leq \bu(t,x)\leq (\phi(x)/\phi(\sigma_1))e^{\rho_0 t}, \hbox{ for all } (t,x)\in[0,\infty)\times\R.
\]
For all $t \geq t_\lambda$ (so that $E_\lambda(t)$ is not empty) and for all $x \in E_\lambda (t),$ there holds
\[
\ds (|x|<\sigma_1) \hbox{ or } \left( |x| \geq \sigma_1 \hbox{ and } \lambda = u(t,x) \leq  (\phi(x)/\phi(\sigma_1) )e^{\rho_0 t}\right).
\]
In all cases, we get that
\[
  \forall t\geq t_\lambda, \forall x\in E_\lambda(t), \phi(x)\geq\min{(\phi(\sigma_1), \lambda \phi(\sigma_1)e^{-\rho_0 t})}.
\]
Then, from the definition of $\phi$ for large $x,$ and since $J(s+\tau)/J(s)\to1$ as $s\to\pm\infty$ for any $\rho>\rho_0/(1-\epsilon_0),$ there exists a time
$\overline{T}_{\lambda} \geq t_\lambda$
such that
\begin{equation}\label{ineq:E_NExp_+}
\forall  t \geq \overline{T}_{\lambda},\  \forall x \in E_\lambda (t), \ph  J(x) \geq  e^{-\rho t},
\end{equation}
which gives~\eqref{eq:ineq_sup_E}.

\

\paragraph*{ \underline{Proof of \eqref{eq:ineq_sup_E} under Hypothesis \ref{hyp:J_geo}} }
Assume that $J$ is satisfies Hyp.~\ref{hyp:J_geo}.
Since $J(x)$ is decreasing with respect to $|x|,$ we get the following inequality for any $x\ge 0:$
\begin{equation}
\begin{array}{rcl}
 \ds\frac{(J\ast J)(x)}{J(x)} & = & \ds\int_{-\infty}^{0}{\frac{J(x-y)}{J(x)}J(y)dy} + \int_{0}^{x}{\frac{J(x-y)J(y)}{J(x)}dy} + \int_{x}^{+\infty}{\frac{J(y)}{J(x)}J(x-y)dy} \\
                         & \leq & \ds1 + \int_{0}^{x/2}{\frac{J(x-y)}{J(x)}J(y)dy} + \int_{x/2}^{x}{\frac{J(x-y)}{J(x)}J(y)dy} \\
                         & \leq & \ds1 + \frac{J(x/2)}{J(x)}.
\end{array}
\label{eq:ineg_J_star}
\end{equation}
From the symmetry of $J,$ the inequality also holds for $x\leq0.$
Moreover, since $J$ satisfies~\eqref{eq:J_geo}, there exists $x_0>0,$ and $C_0>0$ such that
\[
 \ds\ln{\left(\frac{J(x/2)}{J(x)}\right)}=\int_{|x|/2}^{|x|}{\left|\frac{J'(s)}{J(s)}\right|ds}\leq \int_{|x|/2}^{|x|}{\frac{C_0}{s}ds}=C_0\ln{(2)} \hbox{ for all } |x|\geq x_0.
\]
This implies that there exists $K>0$ such that for all $x\in\R$, 
\be \label{eq:J*J/J}
\ds \frac{(J\ast J)(x)}{J(x)}\leq 1+K.
\ee
Since $u_0$ is compactly supported, there exists $\sigma_1>0$ such that $J(\sigma_1)\leq1 $ and $J(x)\geq u_0(x),$ for all $|x|\geq\sigma_1.$
Then, set $\rho_0=r+K$ and for all $(t,x)\in[0,+\infty)\times\R:$
\[
 \bu(t,x)=\min{\left(\frac{J(x)}{J(\sigma_1)}e^{\rho_0 t},1\right)}.
\]
Observe that $u_0(x)\leq\bu(0,x)$ for all $x\in\R.$ Let us now check that $\bu$ is a supersolution of the equation~\eqref{eq:1} satisfied by $u.$ In the region $(t,x)$ such that $\bu(t,x)=1,$ the same arguments 
as
in Sec.~\ref{sec:hyp_J_nonexp} lead to:
\[
 \bu_t(t,x)-J\ast\bu(t,x)+\bu(t,x)-f(\bu(t,x))\geq 1-J\ast\bu(t,x)\geq 0.
\]
Let us check that $\bu$ is also a supersolution of~\eqref{eq:1} when $\bu < 1.$ If $t\geq0,$ $ |x|\geq\sigma_1$ and $\bu(t, x) < 1$ then it follows from~\eqref{eq:J*J/J} that
\begin{equation*}
\begin{array}{rcl}
 \ds (J\ast\bu)(t,x) & \leq & \ds (J\ast J)(x) \frac{e^{\rho_0 t}}{J(\sigma_1)} \\
                   & \leq & \ds  (1+K)\frac{J(x)}{J(\sigma_1)}e^{\rho_0 t} \\
                   & \leq & \ds (1+K)\bu(t,x).
\end{array}
\end{equation*}
This implies that
\begin{equation}
 \begin{array}{rcl}
  \bu_t(t,x)-(J\ast \bu)(t,x)+\bu(t,x)-f(\bu(t,x)) & \geq & \ds \rho_0 \bu(t,x) - (J\ast\bu)(t,x) +\bu(t,x) - r \bu(t,x) \\
                                                   & \geq & \ds (\rho_0 - (1+K) +1 -r) \bu(t,x) \\
                                                   & \geq & 0.
 \end{array}
\end{equation}
The parabolic maximum principle~\cite{Wei82,Yag09} implies that:
\[
\ds u(t,x)\leq \bu(t,x)\leq \frac{J(x)}{J(\sigma_1)}e^{\rho_0 t}, \hbox{ for all } (t,x)\in[0,\infty)\times\R.
\]
For all $t \geq t_\lambda$ (so that $E_\lambda(t)$ is not empty) and all $x \in E_\lambda (t),$ there holds
\[
\ds (|x|<\sigma_1) \hbox{ or } \left( |x| \geq \sigma_1 \hbox{ and } \lambda = u(t,x) \leq \frac{J(x)}{J(\sigma_1)} e^{\rho_0 t}\right).
\]
In all cases, one gets that
\[
  \forall t\geq t_\lambda, \forall x\in E_\lambda(t), J(x)\geq\min{(J(\sigma_1), \lambda J(\sigma_1)e^{-\rho_0 t})},
\]
Then, for any $\rho>\rho_0\geq f'(0)>0,$ there exists a time $\overline{T}_{\lambda} \geq t_\lambda$ such that
\begin{equation}\label{ineq:E_Geo_+}
\forall  t \geq \overline{T}_{\lambda},\  \forall x \in E_\lambda (t), \ph  J (x) \geq  e^{-\rho t},
\end{equation}
which proves~\eqref{eq:ineq_sup_E}.
\carre

\section{Discussion}
We have analyzed the spreading properties of an integro-differential equation with exponentially unbounded or ``fat-tailed'' kernels. Since the pioneering work of Kot et 
al.~\cite{KotLew96}, there have been few mathematical papers on integral equations with exponentially unbounded kernels. However, such slowly decaying kernels are highly relevant in 
the context of population dynamics with long distance dispersal events~\cite{Cla98,ClaFas98}.

We proved that for kernels $J$ which decrease to $0$ slower than any exponentially decaying function the 
level sets of the solution $u$ of the problem~\eqref{eq:1} propagate with an infinite asymptotic speed. This first result shows the qualitative difference between dispersal operators with exponentially unbounded 
kernels and dispersal operators with exponentially bounded kernel which are known to lead to finite spreading speed~\cite{Aro77,Die79,Thi79,Wei82}. This result supports the use of ``fat-tailed'' dispersal 
kernels to model accelerating propagation or fast propagation phenomena~\cite{Cla98,ClaFas98,Ske51}.
 
Moreover, we obtained lower and upper bounds for the position of any level set of $u.$ 
These bounds allowed us to estimate how the solution accelerates, depending on the kernel $J$: the slower the kernel decays, the faster the level sets propagate.
Through several examples, we have seen in Sec.~\ref{sec:main_results} that the level sets of the solution of problem~\eqref{eq:1} move almost linearly when $J$ is close to an exponentially bounded kernel 
(see example~\eqref{exp:NEpLog}, $J(x)= Ce^{-|x|/\ln(|x|)}$ for $|x|\gg1$) while the level sets move exponentially fast when the kernel $J$ has a very fat-tail 
(see example~\eqref{exp:Geo}, $J(x)= C|x|^{-\alpha}$ for $|x|\gg1$). 

It is noteworthy that our results have been derived under assumptions more general than the KPP assumption $f(s)<f'(0)s.$ Indeed, results of Theorems~\ref{theo:infinite_speed},~\ref{theo:speed_estimate_inf} 
and~\ref{theo:speed_estimate_sup} hold with nonlinearities $f$ which may take a weak Allee effect into account: the maximum of the ``per capita growth rate'' $f(s)/s$ is not necessarily reached at $s=0.$ 
In ecological models, the Allee effect can occur for various reasons~\cite{All38}. For instance, at low densities individuals may have trouble finding mates. Our spreading properties in the case with a weak 
Allee effect are in agreement with the numerical results in~\cite{WanKot02} which show that 
exponentially unbounded dispersal kernels can lead to infinite spreading speeds. The conclusion is very different when the nonlinearity $f$ takes a strong Allee effect into accounts, that is if $f(s)<0$ for 
small value of $s.$ Indeed it is proved in this case that problem~\eqref{eq:1} admits traveling wave solutions with constant speed~\cite{BatFif97,Che97,Cov07}. Thus the solutions of problem~\eqref{eq:1} 
with such nonlinearities  have a finite speed of propagation.

One could wonder whether the solution of problem~\eqref{eq:1} with exponentially unbounded kernel converges to some kind of ``accelerated traveling wave solution'', that is a solution $u(t,x)=\phi(x-c(t))$ where
 $c$ is a superlinear function. Our numerical computations suggest that the answer is no (see Fig.~\ref{fig:u_geo_vs_exp}~(a)) since the solutions $u$ becomes flat at large time. The computations also suggest
that the solution $u(t,x)$ is not a generalized transition wave in the sense of~\cite{BerHamP,BerHam07}. 
Such a result has been proved by Hamel and Roques~\cite{HamRoqP2} for the solution of a reaction-diffusion problem with an exponentially unbounded initial condition. In our case, the proof seems to be more 
involved since the operator $u\mapsto J\ast u-u$ is not a differential operator.

 \bibliographystyle{abbrv}

\end{document}